\documentclass[12pt, a4paper, twoside,reqno]{amsart}

\usepackage[margin=2.85cm]{geometry}

\usepackage[T1]{fontenc} 
\usepackage[latin1]{inputenc}
\usepackage{amsmath}
\usepackage{amsthm} 
\usepackage{amssymb}
\usepackage{bbm} 
\usepackage{enumerate} 
\usepackage[dvips]{graphicx}

\usepackage{array}

\theoremstyle{definition}

\theoremstyle{remark}

\theoremstyle{plain}
\newtheorem{thm}{Theorem}[section]
\newtheorem{lem}[thm]{Lemma}
\newtheorem{prop}[thm]{Proposition}

\newtheorem*{khinthm}{Khintchine's Transference Principle}

\newcommand{\norm}[1]{\ensuremath{\left\Vert #1 \right\Vert}}
\newcommand{\abs}[1]{\ensuremath{\left\vert #1 \right\vert}}

\newcommand{\infabs}[1]{\ensuremath{\left\vert #1 \right\vert_\infty}}
\newcommand{\natnum}{\ensuremath{\mathbb{N}}}
\newcommand{\integer}{\ensuremath{\mathbb{Z}}}
\newcommand{\rationals}{\ensuremath{\mathbb{Q}}}

\newcommand{\torus}{\ensuremath{\mathbb{T}}}
\newcommand{\ie}{{\emph{i.e.}}}
\DeclareMathOperator{\dist}{dist}

\DeclareMathOperator{\lcm}{lcm}

\DeclareMathOperator{\skel}{skel}
\DeclareMathOperator{\cosec}{cosec}

\begin{document}

\markboth{M. M. Dodson and S. Kristensen}{Khintchine's theorem and
  transference principle for star bodies}

\title{Khintchine's theorem and transference principle for
star bodies}

\author{M. M. DODSON}

\address{M. M. Dodson, Department of Mathematics, University of York 
  Heslington, York, YO10 5DD, United Kingdom} \email{mmd1@york.ac.uk}

\author{S. KRISTENSEN}

\address{S. Kristensen, School of Mathematics, University of
  Edinburgh, James Clerk Max\-well Building, King's Buildings, Mayfield
  Road, Edinburgh, EH9 3JZ, United Kingdom} 

\curraddr{Department of Mathematical Sciences, Faculty of Science,
  University of Aarhus, Ny Munkegade, Building 530, DK-8000 Aarhus C,
  Denmark}

\email{sik@imf.au.dk}

\keywords{Diophantine approximation; distance functions; metric
  theory.}

\subjclass[2000]{11J83, 11H16}

\begin{abstract}
  Analogues of Khintchine's Theorem in simultaneous Diophantine
  approximation in the plane are proved with the classical height
  replaced by fairly general planar distance functions or equivalently
  star bodies.  Khintchine's transference principle is discussed for
  distance functions and a direct proof for the multiplicative version
  is given.  A transference principle is also established for a
  different distance function.
\end{abstract}

\maketitle
\section{Introduction}
\label{sec:introduction}

A star body $S$ in Euclidean space $ \mathbb{R}^n$ is defined as an open
set containing the origin and for which given any $x\in \mathbb{R}^n$,
there exists a $t_0 \in (0, \infty]$ such that for $t < t_0$, $tx \in
S$ and for $t > t_0$, $tx \notin S$ (see \cite{MR28:1175}). To such
sets, one can associate a continuous distance function $F: \mathbb{R}^n
\rightarrow [0,\infty)$ such that for any $x \in \mathbb{R}^n$ and any $t
\geq 0$, $F(tx) = t F(x)$.  The open star body $S$ may be expressed as
the set of points satisfying $F(x) < 1$. Conversely, we can associate
an open star body to any distance function by this relation.  Some
authors consider only star bodies symmetric about the origin, but this
restriction is not needed here and we consider the more general case.
However, we restrict ourselves principally to the planar case ($n=2$),
as in higher dimensions the number of cases to be considered
proliferates and geometry becomes very complicated.

Examples of planar distance functions include the \emph{height} on
$\mathbb{R}^2$ given by $\infabs{(x_1,x_2)} \allowbreak =
\max\{\abs{x_1}, \abs{x_2}\}$, where the associated star body is a
square. Another example is the function given by $(x_1, x_2) \mapsto
\sqrt{\abs{x_1}\abs{x_2}}$ where the associated star body is bounded
by the hyperbola $x_2 = 1/x_1$ and its mirror image $x_2 = -1/x_1$.
Such symmetric star bodies and their relationship to lattices have
been studied extensively. We refer to \cite{MR28:1175} for an
excellent introduction.

In this paper, we study the relationship between planar star bodies
and metric Diophantine approximation. Let $F : \mathbb{R}^2 \rightarrow
[0, \infty)$ be a distance function and let $\psi: \natnum \rightarrow
(0,\infty)$ be a function with $q \psi(q)$ non-increasing. For
convenience we study the set
\begin{align}
  \label{eq:1}
  W(F;\psi) &= \big\{x \in \torus^2 : F(x-p/q) < \psi(q) \text{ for
    some } p \in \integer^2 \notag \\
  & \qquad \text{ for infinitely many } q \in \natnum\big\} \notag \\
  &= \bigcap_{N=1}^\infty \bigcup_{q=N}^\infty \bigcup_{p \in
    \mathbb{Z}^2} \big\{x \in \torus^2 : F(qx-p) < q\psi(q) \big\},
\end{align}
where $\torus^2$ denotes the two-torus $\{(x_1,x_2) \in \mathbb{R}^2
\colon 0\le x_1,x_2 <1\}$. When $F$ is the height $\infabs{\cdot}$,
$W(F;\psi)$ is the set of simultaneously $\psi$-approximable points.
When $F$ is the function $(x_1, x_2) \mapsto
\sqrt{\abs{x_1}\abs{x_2}}$, $W(F;\psi)$ is the set of multiplicatively
$\psi$-approximable points. In both of these special cases, we can
restrict $p$ to those points with $0 \leq \abs{p} < q$ without loss of
generality. As we shall see, such restrictions are not natural in the
general case.

Let $\abs{E}$ denote the Lebesgue measure of $E \subseteq \mathbb{R}^2$.
Khintchine's Theorem (see \emph{e.g.} \cite{MR50:2084}), which asserts
that when $q \psi(q)$ is decreasing,
\begin{equation}
  \label{eq:4}
  \abs{W(\infabs{\cdot};\psi)} =
  \begin{cases}
    0 & \text{if } \sum_{q=1}^\infty q^2 \psi(q)^2 < \infty, \\
    1 & \text{if } \sum_{q=1}^\infty q^2 \psi(q)^2 = \infty,
  \end{cases}
\end{equation}
was extended by Gallagher~\cite{MR28:1167} who established a general
theorem implying~\eqref{eq:4} and the following:
\begin{equation}
  \label{eq:12}
  \abs{W((x_1, x_2) \mapsto \sqrt{\abs{x_1}\abs{x_2}};\psi)} =
  \begin{cases}
    0 & \text{if } \sum_{q=1}^\infty q^2 \psi(q)^2 \log(
    q^{-1}\psi(q)^{-1}) < \infty,
    \\
    1 & \text{if } \sum_{q=1}^\infty q^2 \psi(q)^2 \log(
    q^{-1}\psi(q)^{-1}) = \infty.
  \end{cases}
\end{equation}
Evidently both these results are of Khintchine type, where the measure
of a set of points $\psi$-approximable by a distance function is
either null or full, according to the convergence or divergence of an
`area' series.

The above results suggest that a general Khintchine-type theorem might
exist for star bodies. While Gallagher's result implies Khintchine
type theorems for a number of star bodies, only convex star bodies and
star bodies of the general shape as the ones in \eqref{eq:12} are
covered. The technical requirement for Gallagher's result (property P)
is that given a point $(p_1, p_2)$ inside the star body $S$, the
entire rectangle with vertices $(0,0)$, $(p_1,0)$, $(p_1, p_2)$ and
$(0, p_2)$ must be a subset of $S$.

It would also seem that a growth condition on the size of the
unbounded parts of the star body would influence the breaking point
(\ie, when the sum converges or diverges).  In an earlier paper
\cite{MR92h:11063}, the first author found a measure zero result of
the above type under a technical covering condition. This work was
subsequently generalised to systems of linear forms
\cite{MR93g:11078}. However, the specific star bodies to which these
results can easily be applied are only minor generalisations of the
examples above.

Some recent results in Diophantine approximation on manifolds due to
Ber\-nik, Kleinbock and Margulis \cite{MR2002g:11102, MR99j:11083}
suggest that there might be a connection between distance functions
and Diophantine approximation and that a
general Khintchine-type theorem might exist for general star bodies.
In Section \ref{sec:khintch-type-theor} below, we will prove
Khintchine type theorems for a large class of distance functions. The
results depend critically on the intrinsic arithmetic properties of
the distance functions as well as their geometric properties. We shall
state the results once we have the notational apparatus in place.

A tantalising result due to Aliev and Gruber~\cite{aliev_gruber} is
in a certain sense dual to the problem studied in this paper and
measures the set of lattices which have points inside a given star
body of infinite volume. They show that given a star body $S \subseteq
\mathbb{R}^n$ of infinite volume, almost all lattices have $n$
linearly independent primitive points inside~$S$.

Associated with the sets $W(\infabs{\cdot};\psi)$ and $W((x_1, x_2)
\mapsto \sqrt{\abs{x_1}\abs{x_2}};\psi)$ are two transference
principles, which relate the simultaneous Diophantine approximation of
the coordinates of $x \in \mathbb{R}^n$ to the linear forms $q \cdot x$,
where $x \in \mathbb{R}^n$ and $q \in \integer^n$. Let $\langle x \rangle
= x-k_x \in [1/2,1/2)^n$, where $k_x \in \mathbb{Z}^n$ denotes the
symmetrised distance from $x \in \mathbb{R}^n$ to the nearest point in the
integer lattice $\integer^n$. The following Transference Principle
(see~\cite{MR80k:10048}; a different formulation is in \cite[Chapter
V, Theorem IV]{MR50:2084}) relates the properties of simultaneous
Diophantine approximation, associated with the height, and those of
the dual inner product form.
\begin{khinthm}
  Let $n \in \natnum$ and let $x \in (\mathbb{R} \setminus
  \rationals)^n$. The following two conditions are equivalent:
  \begin{enumerate}[(i)]
  \item For some $\varepsilon > 0$, there are infinitely many $q \in
    \mathbb{Z}^n$ such that
    \begin{equation*}
      \abs{\langle q \cdot x \rangle} \leq \abs{q}^{-n-\varepsilon}.
    \end{equation*}
  \item For some $\varepsilon' > 0$, there are infinitely many $p \in
  \mathbb{Z}$ such that
  \begin{equation*}
    \norm{qx} \leq \abs{q}^{-(1+\varepsilon')/n}.
  \end{equation*}
  \end{enumerate}
\end{khinthm}
Note that if condition (i) is true for some $\varepsilon > 0$, it is
automatically satisfied for all $\hat{\varepsilon} < \varepsilon$, and
similarly for condition~(ii).

There is an analogous multiplicative version for multiplicative
Diophantine approximation and the inner product form.
Sprind\v{z}uk~\cite[page 69]{MR80k:10048} states it without proof (an
$n$-th is root missing from the left hand side of his inequality~(9)).
Dyson~\cite{dyson47} deduces it from a more general form of
Khintchine's Transference Principle and \cite{WYZ79} gives a
contrapositive form using a method of A.~Baker~\cite{ABaker67}.

\begin{thm}
  \label{thm:multi-trans}
  Let $n \in \natnum$ and $x=(x_1,\dots,x_n) \in (\mathbb{R}\setminus
  \rationals)^n$. The following two conditions are equivalent
  \begin{enumerate}[(i)]
  \item For some $\varepsilon > 0$, there are infinitely
    many $q \in \integer^n$ such that
    \begin{equation}
      \label{eq:15}
      \abs{\langle q \cdot x \rangle} \leq \left(\prod_{i=1}^n
        \max(\abs{q_i},1) \right)^{-1-\varepsilon}.
    \end{equation}
  \item For some $\varepsilon' > 0$, there are infinitely
    many $p \in \integer$ such that
    \begin{equation}
      \label{eq:16}
      \left(\prod_{i=1}^n
        \abs{\langle p x_i \rangle}\right)^{1/n} \leq
      \abs{p}^{-(1+\varepsilon')/n}.
    \end{equation}
  \end{enumerate}
\end{thm}

In Section \ref{sec:transf-princ}, we shall give a direct proof of
Theorem \ref{thm:multi-trans} which also gives a procedure for proving
transference principles for other symmetric distance functions than
$(x_1, \dots ,x_n) \mapsto \left(\prod \abs{x_i}\right)^{1/n}$. We
will also outline an application to a distance function not covered by
previous results. In this case, the assumption that the star bodies
are symmetric about the origin is critical.

We will use the Vinogradov notation. That is, for two real numbers
$x,y$, we will write $x \ll y$ if there exists a constant $C > 0$ such
that $x \leq C y$. If $x \ll y$ and $y \ll x$, we will write $x \asymp
y$.

\section{Khintchine type theorems}
\label{sec:khintch-type-theor}

We will treat two separate classes of distance functions, for which
the theorems are of a different nature. We make the appropriate
definitions and state the results before going on to prove the main
theorems. In the final part of this section, we will discuss the
remaining distance functions and conjectures about the corresponding
Khintchine type theorems.

\subsection{Notation, definitions and statement of results}
\label{sec:notat-defin-stat}
Let $F: \mathbb{R}^2 \rightarrow [0,\infty)$ be a distance function. We define
the \emph{skeleton of $F$} to be the set
\begin{equation}
  \label{eq:8}
  \skel(F) = F^{-1}(0).
\end{equation}
It is easy to see, that the skeleton of a distance function consists
of the origin together with a (possibly infinite and possibly zero)
number of half-lines starting at the origin. For each such line, the
star body has an unbounded component. If there are no such lines, the
star body is bounded. Note that Gallagher's property P is only of
interest when any half-lines in the skeleton lie on one of the
coordinate axes.

Let $L$ be a half-line from the skeleton of $F$. We will call $L$
\emph{significant} if the component around it carries an infinite
amount of mass of the star body, \emph{i.e.}, if for any $M > 0$,
\begin{equation*}
  \left|\left\{x \in \mathbb{R}^2 : F(x) < 1\right\} \cap \left\{x \in
  \mathbb{R}^2 : \dist(x,L) < M \right\} \right|= \infty.
\end{equation*}
Note that the existence of a significant line in the skeleton of the
star body immediately implies that the measure of the star body is
infinite. Conversely, if the measure of an unbounded star body with
only finitely many half-lines in its skeleton is infinite, then the
skeleton contains a significant line.

When the skeleton of $F$ contains an irrationally sloped significant
line, we prove the following result, generalising a result by the
second author \cite{kristensen:_metric_dioph}.
\begin{thm}
  \label{thm:irrational}
  Let $F: \mathbb{R}^2 \rightarrow [0,\infty)$ be a distance function
  corresponding to an unbounded star body and let $\psi: \natnum
  \rightarrow (0,\infty)$ be some function.  Suppose that that at
  least one of the (half-)lines in $\skel(F)$ is significant with an
  irrational slope. Suppose further that for each $\varepsilon > 0$,
  the width $w(r)=w(r;\varepsilon)>0$ of this unbounded component of
  $F^{-1}([0,\varepsilon))$ containing this line is non-increasing as
  the distance $r$ from the origin increases. Then, for almost all $x
  \in \mathbb{R}^2$ and any $q \in \natnum$, there are infinitely many
  $p \in \integer^2$ such that
  \begin{equation}
    \label{eq:3}
    F(x -p) < \psi(q),
  \end{equation}
  whence
  \begin{equation*}
    \abs{W(F;\psi)} = 1.
  \end{equation*}
\end{thm}
Thus the existence of an irrational significant line implies the
stronger conclusion that for almost all $x\in \mathbb{R}^2$ and for
any $\varepsilon$, there are infinitely many $p\in \integer^2$ for
which $F(x-p)<\varepsilon$.  Note that it is possible for a star body
to be unbounded while having no significant lines in the skeleton. We
will return to this case and others like it later. While the
irrationality of the slope makes the result plausible, the result is
not always true and a monotonicity condition is critical in the proof.
The proof is not straightforward and uses recent work of Bugeaud on a
variant of the usual form of Diophantine
approximation~\cite{MR1972699}.

In the case when the skeleton of $F$ consists of finitely many
half-lines whose slopes are rational numbers (or in the degenerate
case when the skeleton is just the origin), we will need additional
definitions. Let $F$ be such a distance function, let $n$ be the
number of lines in the skeleton of $F$ and suppose that the slopes of
these are $s_i/r_i$ for some integers $s_i, r_i$, $i = 1, \dots, n$.
We see that if we take the union of the sets $\{x \in \mathbb{R}^2
: F(x) < \rho\} + p$ where $p \in \integer^2$, we get a pattern which
repeats itself with period $\hat{s} = \lcm(s_1, \dots, s_n)$ in the
direction of the first coordinate axis and with period $\hat{r}=
\lcm(r_1, \dots, r_n)$ in the direction of the second axis.  We call
such a rectangle $\mathcal{R}_F$ a \emph{fundamental rectangle for $F$}.

Let $\varepsilon > 0$ and define the function $D_F$ of $\varepsilon$ by
\begin{equation}
  \label{eq:9}
  D_F(\varepsilon) = \dfrac{\abs{\left(\bigcup_{p \in \integer^2}\{x \in
        \mathbb{R}^2 : F(x) < \varepsilon\} + p \right) \cap
      \mathcal{R}_F}}{\abs{\mathcal{R}_F}}.
\end{equation}
When the star body is bounded, the appropriate definition of the above
quantity is
\begin{equation}
  \label{eq:2}
  D_F(\varepsilon) = \left|\left\{x \in \mathbb{R}^2 : F(x) <
      \varepsilon\right\} \right|.
\end{equation}

We can now state our second Khintchine type theorem.
\begin{thm}
  \label{thm:rational}
  Let $F:\mathbb{R}^2 \rightarrow \mathbb{R}$ be a distance function such that
  $\skel(F)$ consists of finitely many lines, each with a rational
  slope. Let $\psi: \natnum \rightarrow (0,\infty)$ be a function such
  that $D_F(q \psi(q))$ is decreasing. Then
  \begin{equation*}
      \abs{W(F;\psi)} =
  \begin{cases}
    0 & \text{whenever } \sum_{q=1}^\infty D_F(q \psi(q)) < \infty,
    \\
    1 & \text{whenever } \sum_{q=1}^\infty D_F(q \psi(q)) = \infty.
  \end{cases}
  \end{equation*}
\end{thm}
Note that Theorem \ref{thm:rational} covers \eqref{eq:4},
\eqref{eq:12} and a number of additional cases, including all distance
functions corresponding to bounded star bodies as well as all the
cases covered by Gallagher's result \cite{MR28:1167}. Note also, that
$\abs{\mathcal{R}_F}$ depends only on the distance function, and so
the series in the statement of Theorem \ref{thm:rational} could be
replaced by volume sums, as is the case in \cite{MR28:1167}.

\subsection{The irrational case}
\label{sec:irrational-case}

In this section, we prove Theorem \ref{thm:irrational}. Let $q \in
\mathbb{N}$ be fixed.  Let
\begin{equation*}
  L=\{(x,\alpha x)\colon x\in [0,\infty)\}
\end{equation*}
denote a line through the origin in $\skel(F)$ with irrational slope
$\alpha=\tan \theta$.  Suppose that $\alpha > 1$, as the other case
may be treated analogously by interchanging the axes. Identify each
$L$ with its canonical projection into the unit square, which is in
turn identified with the two-torus $\mathbb{T}^2$. Since $L$ has
irrational slope, it follows from Kronecker's Theorem (see
\emph{e.g.}~\cite[Proposition 1.5.1]{MR96c:58055}) that each of these
geodesics is dense in $\torus^2$. However, this is not sufficient to
imply Theorem~\ref{thm:irrational} for significant lines.

Consider a fixed but arbitrary horizontal line $H =\{(x,y_0)\colon
0\le x<1\}$, where $0\le y_0 <1$, through the unit square $\torus^2$;
$H$ will be identified with the circle $\torus^1$ (one dimensional
torus). Consider the geodesic
\begin{equation*}
  L=\left\{\left(\frac{y}{\alpha},y\right)\colon y\in
    [0,\infty)\right\}.
\end{equation*}
The point $((y_0+n)/\alpha, y_0+n)=(y_0/\alpha, y_0)+n(1/\alpha,1)$ on
$L$ is at distance
\begin{equation}
 \label{eq:rn}
  r_n=\left((y_0+n)/\alpha)^2+ (y_0+n)^2\right)^{1/2}
  =(y_0+n)\cosec \theta
\end{equation}
from the origin and projects to the point
\begin{equation*}
  \left(\left\{\frac{y_0+n}{\alpha}\right\},
    \left\{y_0+n\right\}\right)=
  \left(\left\{\frac{y_0+n}{\alpha}\right\},y_0\right)\in \torus^2.
\end{equation*}
As the slope of $L$ is irrational, the geodesic intersects the
horizontal line $H$ first at the point $(x_0,y_0)$, where since
$\alpha>1$,
\begin{equation*}
  x_0=x_0(\alpha)= \dfrac{y_0}{\alpha}<1
\end{equation*}
and thereafter at the distinct points $(T^n x_0,y_0)$, where
$n=1,2,\dots$ and the map $T\colon \torus^1 \rightarrow \torus^1$ can
be regarded as a rotation of the circle by the irrational angle
$\alpha$ and
\begin{equation*}
  T^n x_0:=x_n=\{x_0+n\alpha^{-1}\}.
\end{equation*}

Let $L_n$ be the directed segment which meets
\begin{equation*}
  (T^n x_0,y_0)=(x_n,y_0)=(\{x_0+n\alpha^{-1}\},y_0)
\end{equation*}
from below.  The distance $r_n$ along the irrational line segments
$L_0, L_1, \dots$ from the origin to the point
$(x_n,y_0)=(\{x_0+n\alpha^{-1}\},y_0)$ is given by~\eqref{eq:rn}.

Suppose $L$ a significant line.  Assume for the sake of simplicity
that the star body is \emph{approximately symmetric about significant
  lines}, that is, for any point $(x,y)$ on a significant line $L$,
the distance from $(x,y)$ to the upper part of the boundary
$F^{-1}(\psi(q))$ is comparable to the distance from the lower part.
More precisely, denote by $w^{\pm}=w^{\pm}(r,\psi(q))$ the
perpendicular distances of the point $(x,y)\in L$, where
$r^2=x^2+y^2$, from the upper ($+$) and lower ($-$) boundaries
$\partial F^{-1}([0,\psi(q)))$ of the starbody, so that the width
$w=w(r,\psi(q))$ of the $\psi(q)$-neighbourhood at $(x_n,y_0)$ is
given by $ w=w^+ + w^-$.  Writing $w_n^+=w^+((y_0+n)\cosec
\theta,\psi(q))$ and similarly for $w_n^-$ and $w_n$, we have
\begin{equation*}
  w_n =w_n^+ + w_n^-,
\end{equation*}
where $w_n^+ \asymp w_n^-\asymp w_n$.

\begin{figure}[htbp]
  \centering
  \includegraphics[width=12cm]{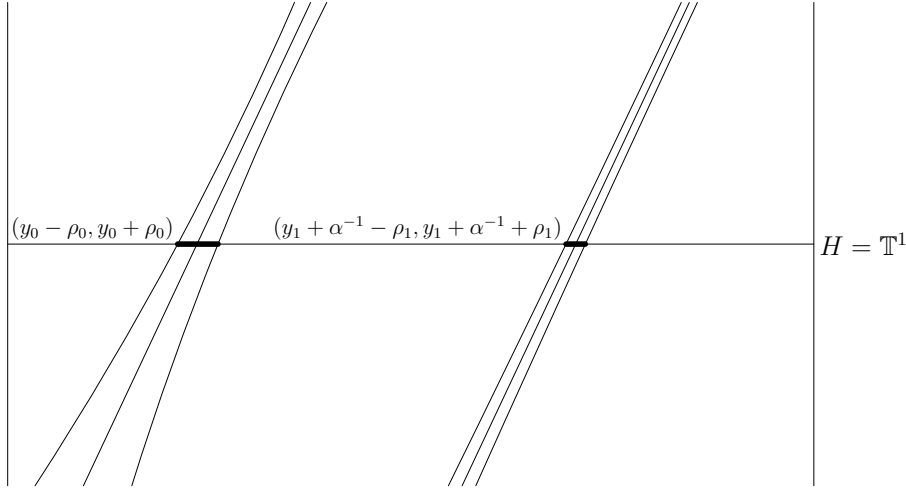}
  \caption{Two consecutive intervals $I_n$ and $I_{n+1}$}
  \label{fig:rotation}
\end{figure}

Let $I_n$ be the interval formed by the intersection of $H$ and the
$\psi(q)$-neigh\-bourhood of the segment $L_n$ and its continuation
through $(x_n,y_0)$ (see Figure~1).  Then the length $|I_n|$ of the
interval $I_n$ satisfies
\begin{equation*}
  w_n\cosec\theta \ll|I_n|\ll w_n\cosec\theta,
\end{equation*}
that is, $w_n\asymp |I_n|$.  Inscribe an interval centred at
$(x_n,y_0)$ with radius
\begin{equation*}
  \sigma_n= K\min\{w^+_n,w^-_n\}.
\end{equation*}
where $K$ is a sufficiently small positive constant chosen so that
\begin{equation}
  \label{eq:rhon}
  \tilde I_n :=\left\{ x \in \mathbb{T}^1 \colon
    \left\Vert x-\left\{x_0-n \alpha^{-1}\right\}\right\Vert
    <\sigma_n\right\}\subseteq I_n,
\end{equation}
where $\Vert x \Vert$ denotes the distance to nearest integer (this is
possible by the assumption of approximate symmetry of the star body).
Then $\sigma_n\asymp w_n$ and $ |\tilde I_n |\asymp |I_n|$.  We show
that given $q \in \mathbb{N}$, almost all $(x,y_0)\in H$ lie in the
interval $\tilde I_n\subseteq I_n$ for infinitely many $n$.

Consider the `limsup' set $E(y_0)$ of points $x$ such that the point
$(x,y_0)$ falls into infinitely many of the intervals $\tilde I_n$ in
$H$.  Now
\begin{equation*}
  \sum_n \abs{\tilde I_n}\asymp \sum_n \abs{ I_n}
\end{equation*}
and
\begin{equation*}
  \sum_n\abs{I_n} \asymp \sum_n \sigma_n  \asymp \sum_n w(r_n,\psi(q))
  \asymp \sum_n w((y+n)\cosec \theta, \psi(q)).
\end{equation*}
But
\begin{equation*}
  \sum_n\abs{I_n}   \asymp \int  w(u\cosec \theta,\psi(q)) du
  \asymp \int  w(v,\psi(q)) dv=\infty,
\end{equation*}
since by assumption $w(u,\psi(q))$ is non-increasing and since the
line $L$ is significant.  Thus
\begin{equation}
  \label{eq:sumdiv}
  \sum_n \abs{\tilde I_n} =\infty.
\end{equation}

We now show that $E(y_0)$ has full measure.
By~\cite{beresnevich03:_measur} it suffices to prove local ubiquity of
the points $x_0+\{n/\alpha\}$, $n=1,2, \dots$ relative to an
appropriate function $\lambda$. This relies on the so-called Three
Distances Theorem due to S\'os \cite{MR20:34} and to a certain extent
follows the methods used by Bugeaud in \cite{MR1972699} (see also
\cite{MR1992153} for an approach using continued fractions).

Let $\lambda:(0,\infty) \rightarrow (0,\infty)$ be a function
decreasing to zero and let $(N_r)$ be a strictly increasing sequence
of natural numbers. The set of points $\mathcal{R} = \{z_n: n \in
\natnum\}$ in the metric space $\mathbb{T}^1$ is said to be
\emph{locally ubiquitous relative to $\lambda$ and $(N_r)$} if there
is a constant $\kappa > 0$ such that for any interval $I=(c-\rho,
c+\rho)$ with $\rho$ sufficiently small and any $r \geq r_0 (I)$,
\begin{equation*}
  \mu\biggl(I \cap \bigcup_{N_r \leq n < N_{r+1}} \left(z_n -
      \lambda(N_r), z_n + \lambda(N_r)\right) \biggr) \geq  2 \kappa
  \rho.
\end{equation*}
This quantitative version of the local density of the points in
$\mathcal{R}$ is the appropriate specialisation of the definition of
local ubiquity in \cite{beresnevich03:_measur}.

\begin{lem}
  \label{lem:ubiquity}
  Let $\mathcal{R} = \{\{x_0 + n \alpha^{-1}\} : n \in \natnum\}$.
  The system $(\mathcal{R}, n)$ is locally ubiquitous relative to the
  function $\lambda(N) = 3/(N_r + 1)$, where $N_{r} \leq N < N_{r+1}$
  for an appropriately chosen strictly increasing sequence of natural
  numbers $(N_r)$.
\end{lem}

\begin{proof}
  As in the proof of \cite[Proposition 1]{MR1972699}, we find that by
  the Three Distances Theorem, for any $N \in \natnum$, the sequence
  of points
  \begin{equation*}
    \{x_0 + \alpha^{-1}\}, \{x_0 + 2 \alpha^{-1}\}, \dots, \{x_0 +
    N\alpha^{-1}\}
  \end{equation*}
  partition the circle $\mathbb{T}^1$ into $N+1$ intervals with
  lengths from a set of cardinality at most three, where $\{x\}$
  denotes the fractional part of $x$. The largest of these distances
  may be seen to be less than or equal to $3/(N+1)$ for infinitely
  many positive integers $N$. Call this sequence $(N_r)$. Hence, for
  any $N \in [N_r, N_{r+1}]$, the intervals centred at the points
  $\{\alpha^{-1}\}, \dots, \{N\alpha^{-1}\}$ of length $3/(N_r+1)$
  cover the circle. Therefore the system is locally ubiquitous
  relative to $\lambda$ and $(N_r)$.
\end{proof}

Consider the family of intervals $\tilde I_n$. The centres of these
intervals form a ubiquitous system relative to $\lambda$ and it is an
easy consequence of~\cite[Theorem 2]{beresnevich03:_measur} that the
divergence of the sum $\sum_n |\tilde I_n|$ implies that the limsup
set $E(y_0)$ has full measure.  It follows that for any fixed $y_0 \in
[0,1)$, for almost there exist infinitely many $p\in\mathbb{Z}^2$ such
that
\begin{equation*}
  F((x,y_0)-p)<\psi(q).
\end{equation*}
Let $\chi_{E(y)}$ denote the characteristic function of $E(y)$.
Theorem \ref{thm:irrational} now follows from Fubini's theorem, since
\begin{equation*}
  \iint_{[0,1)^2} \chi_{E(y)}(x,y) dx dy = \int_0^1 \int_0^1
  \chi_{E(y)}(x,y) dx dy = \int_0^1 1 \ dy = 1,
\end{equation*}
whence
\emph{a fortiori}, for any $q\in \mathbb{N}$,
\begin{equation*}
  \abs{\big\{(x,y) \in \torus^2 : F((x,y)-p) <
    \psi(q) \text{ for infinitely many } p \in \mathbb{Z}^2 \big\}} =
  1.
\end{equation*}

In the more general case, \emph{i.e.}, when the star body is not
assumed to be approximately symmetric about significant lines, the
above argument can be extended by changing the resonant points in the
ubiquitous system to the centres of the intervals $I_n$. These do not
have to be as regularly spaced as the resonant points considered in
the special case, but this difficulty can be resolved by adding an
additional term to the ubiquity function $\lambda$ in Lemma
\ref{lem:ubiquity}. In this case, the ubiquity function should be
\begin{equation*}
  \lambda(N) = 3/(N_r + 1) + \max\{w_{N_r}^+, w_{N_r}^-\}.
\end{equation*}
If necessary, we may take a smaller star body inside the original one
for which $\max\{w_{N_r}^+, w_{N_r}^-\}$ is not dominating the main
term, but for which the corresponding line remains significant. The
details are left to the interested reader.

Note that if the unbounded starbody has no significant lines (so that
$F^{-1}([0,\varepsilon))$ has finite volume), then given $q\in
\mathbb{N}$, the set of points $(x,y)\in \torus^2$ for which
\begin{equation*}
  F((x,y)-p)<\psi(q)
\end{equation*}
for infinitely many $p\in\mathbb{Z}^2$ is null by the Borel--Cantelli
lemma.  The question of the measure of the set of $(x,y)$ for which
\begin{equation*}
  F\left((x,y)-\dfrac{p}{q}\right)<\psi(q)
\end{equation*}
holds for infinitely many $p/q$ is open.

\subsection{The rational case}
\label{sec:rational-case}

For completeness, we first deal with the degenerate case when the
skeleton consists solely of the origin, corresponding to $F$ being a
gauge function~\cite{MR28:1175}.
\begin{lem}
  \label{lem:bounded}
  Let $S \subseteq \mathbb{R}^2$ be a bounded star body with corresponding
  distance function $F$. Then
  \begin{equation*}
      \abs{W(F;\psi)} =
  \begin{cases}
    0 & \text{if } \sum_{q=1}^\infty D_F(q \psi(q)) < \infty,
    \\
    1 & \text{if } \sum_{q=1}^\infty D_F(q \psi(q)) = \infty.
  \end{cases}
  \end{equation*}
\end{lem}

\begin{proof}
  By the final corollary of \cite[Section IV.2]{MR28:1175}, there
  are constants $c, C > 0$ such that for any $x \in \mathbb{R}^2$,
  \begin{equation}
    \label{eq:10}
    c \infabs{x} \leq F(x) \leq C \infabs{x},
  \end{equation}
  whence
  \begin{equation*}
    W\left(\infabs{\cdot};\dfrac{\psi}{C}\right) \subseteq W(F;\psi)
    \subseteq W\left(\infabs{\cdot};\dfrac{\psi}{c}\right).
  \end{equation*}
    Applying \eqref{eq:4} immediately gives the result.
\end{proof}

We now suppose that the skeleton of the star body contains a half-line
and prove the two cases of Theorem~\ref{thm:rational} separately.
First, for any $q \in \natnum$ and any $\varepsilon > 0$, we define
the resonant neighbourhood
\begin{equation}
  \label{eq:17}
  B_q(F, \varepsilon) = \bigcup_{p \in \integer^2}\left\{ x \in \torus^2
    : F (x -p/q) < \varepsilon \right\}.
\end{equation}
The measure of the resonant neighbourhood is closely related to
$D_F(\varepsilon)$. Note first that
\begin{equation}
  \label{eq:14}
  \abs{\tfrac{1}{q} \mathcal{R}_F \cap B_q(F, \varepsilon)} =
  \dfrac{\hat{s} \hat{r}}{q^2} D_F(q \varepsilon),
\end{equation}
as is seen by scaling the set to be estimated by $q$ and using
homogeneity of the distance function.

To obtain the measure of $B_q(F, \varepsilon)$, we tile $\mathbb{R}^2$
by disjoint translates of $\tfrac{1}{q}\mathcal{R}_F$. Inside each of
these disjoint sets, we find a translated copy of $\tfrac{1}{q}
\mathcal{R}_F \cap B_q(F, \varepsilon)$, and the union of these sets
cover $B_q(F, \varepsilon)$. Hence, to estimate the measure of $B_q(F,
\varepsilon)$, it suffices to count the maximal (resp. minimal) number
of disjoint translates of $\tfrac{1}{q}\mathcal{R}_F$ that can fit
inside (resp. are needed to cover) the unit square, and multiply the
result by the estimate from \eqref{eq:14}. In this way, we get
\begin{equation*}
  \left[\dfrac{q}{\hat{r}}\right]  \left[\dfrac{q}{\hat{s}}\right]
  \dfrac{\hat{s} \hat{r}}{q^2} D_F(q \varepsilon) \leq \abs{B_q(F,
    \varepsilon)} \leq \left(\dfrac{q}{\hat{s}} +1\right) 
  \left(\dfrac{q}{\hat{r}} +1\right) \dfrac{\hat{s}\ \hat{r}}{q^2}
  D_F(q \varepsilon),
\end{equation*}
whence
\begin{equation}
  \label{eq:5}
  \abs{B_q(F, \varepsilon)} \asymp D_F(q \varepsilon).
\end{equation}

Applying the Borel--Cantelli Lemma \cite[Lemma 3.14]{MR93d:60001} to
the resonant neighbourhoods with the measure estimate \eqref{eq:5}
yields a condition for $\abs{W(F;\psi)}$ to be zero. Indeed, whenever
\begin{equation}
  \label{eq:6}
  \sum_{q = 1}^\infty D_F(q \psi(q)) < \infty,
\end{equation}
$W(F;\psi)$ must be null. We have shown:
\begin{lem}
  \label{lem:convergence}
  Let $\psi: \integer \rightarrow (0,\infty)$ be some decreasing
  function and let $F: \mathbb{R}^2 \rightarrow \mathbb{R}$ be a distance
  function such that the lines in $\skel(F)$ all have rational slopes.
  Let $D_F(\varepsilon)$ be defined as above. Suppose that
  \begin{equation*}
    \sum_{q = 1}^\infty D_F(q \psi(q)) < \infty.
  \end{equation*}
  Then $\abs{W(F;\psi)} = 0$.
\end{lem}

We now aim to obtain the corresponding divergence result, \emph{i.e.},
when the series \eqref{eq:6} diverges, we expect the measure of $W(F;
\psi)$ to be full. Suppose first that the central part of the star
body carries the bulk of the mass, \emph{i.e.}, there is a distance
function $\widetilde{F}$ such that
\begin{enumerate}[(i)]
\item  $\widetilde{F}$ determines a bounded star body. \label{item:5}
\item $F(x) \leq \widetilde{F}(x)$ for all $x$ (or equivalently,
  the star body defined by $\widetilde{F}$ is a subset of the one
  defined by $F$ by \cite[Corollary, p. 107]{MR28:1175}). \label{item:6}
\item $D_{\widetilde{F}}(\varepsilon) \asymp D_F(\varepsilon)$ for
  $\varepsilon > 0$ small enough. \label{item:7}
\end{enumerate}
In this case, $W(\widetilde{F};\psi) \subseteq W(F;\psi)$, so the
conclusion of Theorem \ref{thm:rational} is ensured by Lemma
\ref{lem:bounded}. We have shown

\begin{lem}
  \label{lem:bulk_at_centre}
  Let $F$ be a distance function and suppose that $\skel(F)$ consists
  of lines with rational slopes. Suppose further that there is a
  distance function $\widetilde{F}$ satisfying
  (i)--(iii) above.  Then Theorem
  \ref{thm:rational} holds for this distance function.
\end{lem}

Suppose now that the central part of the star body does not carry the
bulk of the mass, so that there is a cusp which includes an infinite
amount of the mass. We will construct a subset of $W(F;\psi)$, which
is better behaved but which still has full measure. In brief, we will
truncate the star body in such a way that only the unbounded component
which carries most of the mass will remain. Subsequently, we will take
the limsup set of these truncated sets, for which we can calculate the
measure.  For convenience we denote by $F^*_\varepsilon$ the distance
function associated with the star body $\{\dist(x, L) < w(\varepsilon)\}
\cap \{F(x) < \varepsilon\}$.

\begin{lem}
  \label{lem:truncation}
  Let $F$ be a distance function as in the statement of Theorem
  \ref{thm:rational}. Then there exists a function $w: (0,\infty)
  \rightarrow (0,\infty)$ with $w(\varepsilon)$ tending to zero as
  $\varepsilon$ tends to zero together with a (half-)line $L \in
  \skel(F)$ such that
  \begin{equation}
    \label{eq:7}
    D_{F^*_\varepsilon}(1) \asymp D_F(\varepsilon),
  \end{equation}
  and so that the intersection between any line perpendicular to $L$
  and the star body associated to $F^*_\varepsilon$ is a line segment.
\end{lem}

\emph{Remark:} We have to choose an entire family of distance
functions rather than just a single one. This is to preserve the
scaling properties of the original star body. This takes us away from
the problem of approximation with respect to distance functions, but
as we shall see, only the shape of the sets studied are of importance
for the remainder of the proof.

\begin{proof}
  We choose the line in $\skel(F)$ about which a significant
  proportion of the mass is concentrated. Such a line must exist,
  since there are only finitely many lines in $\skel(F)$. By
  truncating the star body about this line, it is now straightforward
  to choose the function $w$ in such a way that~\eqref{eq:7} holds. To
  ensure that the second property holds, we choose $w$ as small as
  possible. For at least one of the lines in the $\skel(F)$, we will
  have both properties satisfied. Indeed, otherwise most of the mass
  would be concentrated in the cental part of the star body associated
  to $F$.
\end{proof}

\begin{lem}
  \label{lem:decreasing}
  Let $F'_1, F'_2: \mathbb{R}^2 \rightarrow \mathbb{R}$ be distance
  functions with $\skel(F'_1) = \skel(F_2') = L$, a (half-)line. Let
  $\delta_1, \delta_2 > 0$. Let $v \in \mathbb{R}^2$ be a unit vector
  in the direction of $L$ and let $v^\perp$ be a unit normal to $L$.
  Let $t_1, t_2 \in \mathbb{R}$. The measure
  \begin{equation*}
    \abs{\{x \in \mathcal{R}_{F_1'} : F_1'(x) < \delta_1\} \cap \{x
      \in \mathcal{R}_{F_2'} : F_2'(x+t_1 v + t_2 v^\perp)<\delta_2\}}
  \end{equation*}
  decreases as $\abs{t_i}$ increases and the other coordinate is
  fixed.
\end{lem}

\begin{proof}
  We define the characteristic functions $\chi_1$ and $\chi_2$ of $\{x
  \in \mathcal{R}_{F_1'} : F_1'(x) < \delta_1\}$ and $\{x \in
  \mathcal{R}_{F_2'} : F_2'(x) < \delta_2\}$ respectively. Clearly,
  for fixed $t_1$, the function $\chi_1(x) \chi_2(x+t_1 v + t_2
  v^\perp)$ is the characteristic function of an interval in the $t_2$
  coordinate, and so decreasing in $\abs{t_2}$ and analogously for
  fixed $t_2$. Hence,
  \begin{multline*}
    \abs{\{x \in \mathcal{R}_{F_1'} : F_1'(x) < \delta_1\} \cap \{x
      \in \mathcal{R}_{F_2'} : F_2'(x+t_1 v + t_2 v^\perp)<\delta_2\}}
    \\
    = \int_{\mathbb{R}^2} \chi_1(x) \chi_2(x+t_1 v + t_2 v^\perp) dx
  \end{multline*}
  also decreases with $\abs{t_i}$ when the other coordinate is fixed.
\end{proof}

We will now estimate the overlap between two resonant neighbourhoods
in terms of the measure of the star bodies defined in Lemma
\ref{lem:truncation}. To ease notation, let $B^*_q$ be defined as
$B_q(F,\psi(q))$ in \eqref{eq:17} under the additional restriction
that $(p_1 \hat{r}, q) = ( p_2\hat{s},q) = 1$. Furthermore, let
\begin{equation*}
  A^*_q = \left\{x \in \mathbb{T}^2 : F_{q\psi(q)}^*(x) < 1 \right\},
\end{equation*}
where $F_{q\psi(q)}^*$ is the distance function from Lemma
\ref{lem:truncation}.

\begin{lem}
\label{lem:quasi-independent}
  Let $q, q' \in \mathbb{N}$. Then
  \begin{equation*}
    \abs{B^*_q \cap B^*_{q'}} \leq\abs{A^*_q}\abs{A^*_{q'}}  \asymp
    D_F(q\psi(q)) D_F(q' \psi(q')).
  \end{equation*}
\end{lem}

\begin{proof}
  The last asymptotic equality follows at once from Lemma
  \ref{lem:truncation}. We concentrate on the first inequality.
  First, we note that
  \begin{equation*}
    B^*_q \subseteq \bigcup \dfrac{A_q + \binom{k_1
        \hat{r}}{k_2 \hat{s}}}{q} ,
  \end{equation*}
  where the union is taken over all integers $k_1, k_2$  with
  $0 \leq k_1 \hat{r} \leq q$ and $0 \leq l_1 \hat{r}$
  such that $(q, k_1 \hat{r}) = (q, k_2 \hat{s}) = 1$.
 Hence
  \begin{equation*}
    B^*_q \cap B^*_{q'} \subseteq \bigcup \left(\dfrac{A_q +
        \binom{k_1 \hat{r}}{k_2 \hat{s}}}{q} \cap \dfrac{A_{q'} +
        \binom{l_1 \hat{r}}{l_2 \hat{s}}}{q'}\right),
  \end{equation*}
  where the union is taken over all integers $k_1, k_2, l_1, l_2$ with
  $0 \leq k_1 \hat{r}, k_2 \hat{s} \leq q$ and $0 \leq l_1 \hat{r},
  l_2 \hat{s} \leq q'$ such that $(q, k_1 \hat{r}) = (q, k_2 \hat{s})
  = (q', l_1 \hat{r}) = (q', l_2 \hat{s}) = 1$. Consequently,
  \begin{equation}
    \label{eq:11}
    \begin{split}
      \abs{B^*_q \cap B^*_{q'}} &\leq \sum \abs{\dfrac{A_q + \binom{k_1
            \hat{r}}{k_2 \hat{s}}}{q} \cap \dfrac{A_{q'} + \binom{l_1
            \hat{r}}{l_2
            \hat{s}}}{q'}}\\
      &= \sum \abs{\dfrac{A_q}{q} \cap \left(\dfrac{A_{q'}}{q'} +
          \dfrac{1}{qq'} \binom{(l_1 q - k_1 q')\hat{r}}{(l_2 q - k_2
            q')\hat{s}}\right)},
    \end{split}
  \end{equation}
  where the summation range is as above.  Under the summation
  conditions, there are at most $(q, q')^2$ repeated summands in the
  sums in \eqref{eq:11}. Hence, on changing the summation ranges, we
  get
  \begin{equation}
    \label{eq:13}
    \abs{B^*_q \cap B^*_{q'}} \leq\sum_{(j_1, j_2) \in \mathbb{Z}^2}
    \abs{\dfrac{A_q}{q} \cap \left(\dfrac{A_{q'}}{q'} + \dfrac{(q,
          q')}{qq'} \binom{j_1 \hat{r}}{j_2\hat{s}}\right)}.
  \end{equation}

  Consider now the function
  \begin{equation*}
    f(t_1, t_2) = \abs{\dfrac{A_q}{q} \cap \left(\dfrac{A_{q'}}{q'} +
        \dfrac{(q, q')}{qq'} \binom{t_1 \hat{r}}{t_2\hat{s}} \right)}.
  \end{equation*}
  As a consequence of Lemma \ref{lem:decreasing}, this function
  decreases as $\abs{t_i}$ increases ($i=1,2$). Conseqently,
  \begin{equation*}
    \int_0^1 \int_0^1 f(t_1, t_2) dt_1 dt_2 \geq f(1,1).
  \end{equation*}
  Using this and translation invariance of Lebesgue measure in
  (\ref{eq:13}), we obtain
  \begin{alignat*}{2}
    \abs{B^*_q \cap B^*_{q'}} &\leq (q, q')^2 \sum_{(j_1, j_2) \in
      \mathbb{Z}^2} \iint_{[0,1)^2 + (j_1, j_2)} f(t_1, t_2) dt_1
    dt_2\\
    &= (q, q')^2 \iint_{\mathbb{R}^2} \abs{\dfrac{A_q}{q} \cap
      \left(\dfrac{A_{q'}}{q'} + \dfrac{(q, q')}{qq'} \binom{t_1
          \hat{r}}{t_2\hat{s}} \right)} dt_1 dt_2 \\
    &= \iint_{A_q} \abs{A_{q'}} dt_1 dt_2\\
    &= \abs{A_q} \abs{A_{q'}}.
  \end{alignat*}
  This completes the proof.
\end{proof}

The following estimate uses elementary number theory and summation by
parts (see \cite[Chapter VII, Lemma 7]{MR50:2084}).
\begin{lem}
  \label{lem:euler-phi}
  Let $\omega(q)$ be monotonely decreasing and positive. Let $\phi(q)$
  denote the Euler $\phi$-function of $q$. Then
  \begin{equation*}
    \sum_{q=1}^N \left(\dfrac{\phi(q)}{q}\right)^2 \omega(q) \gg \sum_{q=2}^N
    \omega(q).
  \end{equation*}
\end{lem}

The above lemma will be useful because of the following:
\begin{lem}
  \label{lem:restricted}
  Suppose that $q \in \natnum$ satisfies $(q,\hat{r}) = (q,\hat{s}) =
  1$. Then,
  \begin{equation*}
    \abs{B^*_q} \gg \left(\dfrac{\phi(q)}{q}\right)^2 D_F(q\psi(q)).
  \end{equation*}
\end{lem}

\begin{proof}
  This follows from the definition of $B^*_q$. Indeed, the number of
  sets under the union in the definition is $\asymp (\phi(q)/q)^2$, as
  is easily seen by counting the number of $(p_1, p_2)$ satisfying the
  co-primality condition. The co-primality condition and the fact that
  $D_F(q \psi(q))$ is non-increasing implies that the disjoint 
  components are dominant.
\end{proof}

Note that Lemma \ref{lem:euler-phi} and Lemma \ref{lem:restricted}
together with the divergence assumption imply that
\begin{equation}
  \label{eq:18}
  \sum_{\substack{q\geq 1 \\ (q,\hat{r}) = (q,\hat{s}) = 1}}
  \abs{B^*_q} = \infty.
\end{equation}
Indeed, the set of $q \in \mathbb{N}$ with $(q,\hat{r}) = (q,\hat{s})
= 1$ contains the arithmetic progression
$\{\hat{r}\hat{s}n+1\}_{n=0}^\infty$, and so has positive density. On
noting that the terms of the series are decreasing over the full range
of $q$, \eqref{eq:18} follows by a simple change of variables.

\begin{lem}
  \label{lem:divergence}
  Let $F: \mathbb{R}^2 \rightarrow [0,\infty)$ be a distance function such
  that the lines in $\skel(F)$ have rational slopes. Let $\psi:
  \natnum \rightarrow  [0,\infty)$ be a function such that $D_F(q
  \psi(q))$ is non-increasing and suppose that
  \begin{equation*}
    \sum_{q = 1}^\infty D_F(q\psi(q)) = \infty.
  \end{equation*}
  Then $\abs{W(F;\psi)} = 1$.
\end{lem}

\begin{proof}
  Denote by $\sum_{q,q'}^N$ the sum over all $q, q'\in \mathbb{Z}$
  such that $1\le q, q'\le N$, $(q,\hat{r})=(q,\hat{s})=1$,
  $(q',\hat{r})=(q',\hat{s})=1$, and similarly by $\sum_{q}^N$ the sum
  over all $q \in \mathbb{Z}$ such that $1\le q \le N$,
  $(q,\hat{r})=(q,\hat{s})=1$.  Applying a quasi-independent
  divergence case of the Borel--Cantelli Lem\-ma (see Lemma 5 in
  \cite{MR80k:10048}), we see that
  \begin{alignat*}{2}
    \abs{W(F;\psi)} &\geq \limsup_{N \rightarrow \infty}
    \dfrac{\left(\sum_q^N \abs{B^*_q}\right)^2}{\sum_{q, q'}^N
      \abs{B^*_q \cap B^*_{q'}}}\\
    &\gg \limsup_{N \rightarrow \infty} \dfrac{\left(\sum_q^N
        \left(\dfrac{\phi(q)}{q}\right)^2
        \abs{B_q(F,\psi(q))}\right)^2}{\sum_{q, q'}^N \abs{A_q}
      \abs{A_{q'}}}\\
    &\gg \limsup_{N \rightarrow \infty} \dfrac{\left(\sum_q^N
        D_F(q \psi(q))\right)^2}{\left(\sum_q^N D_F(q
          \psi(q))\right)^2} \gg 1
  \end{alignat*}
  by Lemmas~\ref{lem:quasi-independent},~\ref{lem:restricted}
  and~\ref{lem:euler-phi}. As $W^*(F;\psi) \subseteq W(F;\psi)$, we
  clearly have $\abs{W(F;\psi)} \geq c > 0$ for some $c > 0$.
  
  We now apply an `inflation' argument taken from~\cite[Chapter
  VII]{MR50:2084}.  Let $\eta: \mathbb{N} \rightarrow (0,1]$ be a
  function which decreases to zero so slowly that $\sum_q D_F(q
  \eta(q) \psi(q)) = \infty$.  Such a function can be found as
  follows. First, define a strictly inreasing sequence $q_j$ such that
  $\sum_{q_j\le q< q_{j+1}} D_F(q\psi(q))>1$ for any $j \in
  \mathbb{N}$.  Subsequently, define a sequence of positive real
  numbers $\eta_j$ by the equation
  \begin{equation*}
    D_F(\eta_j q_{j+1} \psi(q_{j+1})) = \tfrac{1}{j} D_F(q_{j+1}
    \psi(q_{j+1})).
  \end{equation*}
  Such numbers exist by continuity of the distance function and since $q
  \psi(q)$ is non-increasing. Finally, let
  \begin{equation*}
    \eta(q)=\eta_j, \ q_j \leq q<q_{j+1}.
  \end{equation*}
  In is now immediate that
  \begin{align*}
    \sum_{q=1}^\infty D_F(q \eta(q) \psi(q)) &= \sum_{j=1}^\infty
    \sum_{q_j\le q< q_{j+1}} D_F(q \eta(q) \psi(q))\\ 
    &\geq \sum_{j=1}^\infty \sum_{q_j\le q< q_{j+1}} D_F(\eta_j
    q_{j+1} \psi(q_{j+1})) > \sum_{j=1}^\infty \dfrac{1}{j} = \infty.
  \end{align*}
  
  Let the function $\psi_\eta: \mathbb{N} \rightarrow [0,\infty)$ be
  defined by $\psi_\eta (q) = \eta(q) \psi(q)$. Going through the
  above argument with $\psi_\eta$ in place of $\psi$, we easily find
  that $\abs{W(F;\psi_\eta)} > c > 0$. Let $x_0 \in W(F;\psi_\eta)$ be
  a point of metric density for $W(F;\psi_\eta)$ and let $\varepsilon
  \in (0,1)$ be arbitrary. By Lebesgue's Density Theorem, there is an
  $n_0 \in \natnum$ and a box $H$ centred at $x_0$ of measure $\abs{H}
  = 1/n_0$ such that
  \begin{equation*}
    \dfrac{\abs{W(F;\psi_\eta)\cap H}}{\abs{H}} =
    n_0\abs{W(F;\psi_\eta)\cap H} \geq 1 - \varepsilon.
  \end{equation*}
  Thus $\abs{n_0(W(F;\psi_\eta)\cap H)} \geq 1 - \varepsilon$, so
  there is a set $U \in \torus^2$ of measure at least $1 -
  \varepsilon$ such that every point $x \in U$ is of the form $x = n_0
  x' + p$, $p \in \integer^2$, $x' \in W(F;\psi_\eta)$.  We claim that
  $U\subseteq W(F;n_0 \psi_\eta)$. Indeed, let $x \in U$. Then
  \begin{equation*}
    F\left(x - \dfrac{p'}{q}\right) = F\left(n_0x'+p-\dfrac{p'}{q}\right) 
    = n_0 F\left(x' - \dfrac{qp-p'}{n_0 q}\right).
  \end{equation*}
  As $p' \in \integer^2$ varies freely, ratio $(p'- pq)/n_0 q$ varies
  freely over the lattice $1/n_0 q \integer^2$ which contains $1/ q
  \integer^2$. Hence, as $x' \in W(F;\psi_\eta)$, for infinitely many
  $q \in \integer$ there is a $p' \in \integer^2$ such that
  \begin{equation*}
    F(x - p'/q) < n_0 \eta(q) \psi(q).
  \end{equation*}
  But there exists $Q$ such that $q\geq Q$ implies that $\eta(q)\le
  n_0^{-1}$, whence
  \begin{equation*}
    F(x- p'/q) < \psi(q)
  \end{equation*}
  for infinitely many $q$ and so $x\in W(F;\psi)$. Hence 
  \begin{equation*}
    \abs{W(F; \psi)} \geq \abs{U} \geq  1-\varepsilon
  \end{equation*}
  and since $\varepsilon$ was arbitrary, we are done.
\end{proof}

Combining Lemmas~\ref{lem:bounded},~\ref{lem:convergence}
and~\ref{lem:divergence} proves Theorem~\ref{thm:rational}.

\subsection{Missing cases}
\label{sec:missing-cases}

Despite  having dealt with a very large class of
distance functions in the preceding two sections, there are still two
cases not covered by our results. The first  is when $\skel(F)$
consists of infinitely many rationally sloped lines, when the notion
of a fundamental region makes little sense. We conjecture that this
case may be treated analogously with the irrational case. This should
be possible as the resonant sets will induce geodesics on the torus
winding around arbitrarily many times.

The second is when each irrationally sloped lines of $\skel(F)$ is
insignificant. We conjecture that these may be removed from the star
body without affecting the corresponding Khintchine type result, which
reduces to the rational case.  While these conjectures seem reasonable
to us, we have not yet been able to prove them.

\section{Transference principles}
\label{sec:transf-princ}

In the rational case, it makes sense to talk about transference
principles. We will give a direct proof of
Theorem~\ref{thm:multi-trans} and then give a description of how this
method may be extended to other star bodies where the skeleton
consists of lines with rational slopes. We prove the theorem in
arbitrary dimension, as this causes no loss of clarity. In addition to
the proof of Theorem \ref{thm:multi-trans}, we illustrate the
versatility of our method by sketching a proof of a transference
principle for a distance function not covered by any previous results.

We first prove some preliminary equivalences.

\begin{prop}
  \label{prop:mult2boxes}
  Let the distance function $F\colon \mathbb{R}^n\to [0,\infty)$ be given by
  \begin{equation*}
    F(x)= F(x_1,\dots, x_n) = \left(\prod_{i=1}^n
      \abs{x_i}\right)^{1/n}
  \end{equation*}
  and define the distance function $H_{\nu} \colon \mathbb{R}^n\to
  [0,\infty)$ by
  \begin{equation}
    \label{oldeq:6}
    H_{\nu}(x) = \max\left\{\nu_1 \abs{x_1}, \dots, \nu_n \abs{x_n}\right\},
  \end{equation}
  where $\nu=(\nu_1, \dots, \nu_n)\in (0,\infty)^n$. Then for each
  $\lambda > 0$, $F(x) \leq \lambda$ if and only if $H_{\nu}(x) \leq
  \lambda$ for some $\nu$ with $\nu_1 \cdots \nu_n = 1$.
\end{prop}

\begin{figure}[htbp]
  \centering
  \includegraphics[width=9cm]{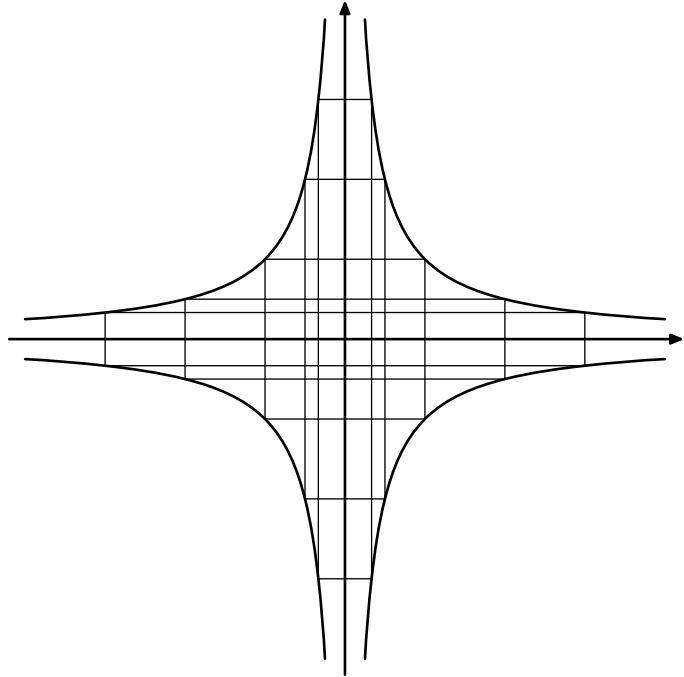}
  \caption{The star body corresponding to $F$ for $n=2$ with some
    parametrising boxes.}
  \label{fig:multiplicative}
\end{figure}

Both functions in Proposition \ref{prop:mult2boxes} are distance
functions which are symmetric about the origin (a planar configuration
is shown in Figure~\ref{fig:multiplicative}).

\begin{proof}
  To begin, assume that $\abs{x_1} \cdots \abs{x_n} \leq \lambda^n$.
  We will find $\nu_1, \dots, \nu_n$ with the required properties.
  For $i = 1, \dots, n-1$, let $\nu_i = \lambda / \abs{x_i}$ and let
  $\nu_n = 1 /(\nu_1 \cdots \nu_{n-1})$. We only need to prove that
  $\nu_n \abs{x_n} \leq \lambda$. But this is the case, since by
  assumption
  \begin{equation*}
    \nu_n \abs{x_n} = \dfrac{\abs{x_1} \cdots
      \abs{x_n}}{\lambda^{n-1}} \leq \dfrac{\lambda^n}{\lambda^{n-1}}
    = \lambda. 
  \end{equation*}
  Then converse is even easier, since
  \begin{equation*}
    \abs{x_1} \cdots \abs{x_n} = \nu_1 \abs{x_1} \cdots \nu_n
    \abs{x_n} \leq \lambda^n
  \end{equation*}
  by assumption.
\end{proof}

\begin{prop}
  \label{prop:pre_transfer}
  Let $\mu, \lambda > 0$. The following are equivalent:
  \begin{enumerate}[(i)]
  \item There is an integer solution $q \in \mathbb{Z}^n
    \setminus \{0\}$ to
    \begin{equation}
      \label{eq:item3}
      \abs{\langle q \cdot x \rangle} \leq \lambda, \quad
      \left(\prod_{i=1}^n \max(\abs{q_i},1)\right)^{1/n} \leq \mu.
    \end{equation}
  \item  There is an integer solution $p \in \mathbb{Z}
    \setminus \{0\}$ to
    \begin{equation}
      \label{eq:item4}
      \left(\prod_{i=1}^n \abs{\langle p x_i\rangle}\right)^{1/n} \leq
      n \lambda, \quad \abs{p} \leq n \mu \lambda^{(1-n)/n}.
  \end{equation}
  \end{enumerate}
\end{prop}

\begin{proof}
  Using Proposition \ref{prop:mult2boxes}, we may encode solutions to
  these equations in matrices as follows. Define the matrices
  \begin{align}
    \label{oldeq:7}
    A &=\left(
      \setlength\extrarowheight{5pt}
    \begin{array}{cccc}
      \lambda^{-1} x_1 & \nu_1 \mu^{-1} & \dots & 0 \\
      \vdots &  \vdots & \ddots & \vdots \\
      \lambda^{-1} x_n & 0  & \dots  & \nu_n \mu^{-1} 
      \rule[-8pt]{0pt}{5pt}
      \\
      \lambda^{-1} & 0 & \cdots & 0
    \end{array}\right), \notag \\
    A^* &=\left( 
     \setlength\extrarowheight{5pt}
    \begin{array}{cccc}
      \lambda & -\mu \nu^{-1}_1 x_1 & \cdots & - \mu \nu^{-1}_n x_n
      \rule[-8pt]{0pt}{5pt}
      \\
     0 & \mu \nu^{-1}_1 &\dots  & 0 \\
     \vdots &\vdots  & \ddots &\vdots \\
      0 & 0 &\cdots & \mu \nu^{-1}_n
    \end{array}\right).
  \end{align}
  If (\ref{eq:item3}) has a non-trivial integer solution $q =
  (q_1,\dots, q_n)$, then for some choice of $\nu_1, \dots, \nu_n$
  with $\nu_1 \cdots \nu_n = 1$, then the inequality
  \begin{equation}
    \label{eq:qA}
  \abs{\tilde{q} A}_\infty \leq 1
  \end{equation}
  has a non-trivial integer solution $\tilde{q} =
  (q_1, \dots, q_n, p)$. This follows on considering each coordinate
  of $\tilde{q} A$ and applying Proposition \ref{prop:mult2boxes}.
  
  As can be verified by calculation, when $(\tilde a,\tilde b) \in
  \mathbb{Z}^{n+1} \times \mathbb{Z}^{n+1}, $ the linear form
  \begin{equation*}
    \label{eq:image}
    \Phi(\tilde a,\tilde b) = \tilde a A  \cdot \tilde b A^*
    =  a_1 b_2 + \dots + a_n b_{n+1} + a_{n+1} b_1\in \integer 
  \end{equation*}
  when $\tilde x,\tilde y$ are integer vectors.  By~\eqref{eq:qA}
  and~\cite[Chapter V, Theorem I]{MR50:2084}, there is an integer
  solution to
  \begin{equation*}
    \abs{\tilde p A^*} \leq n \abs{\det(A^*)}^{1/n} = n \mu
    \lambda^{1/n}, 
  \end{equation*}
  since $\nu_1^{-1} \cdots \nu_n^{-1} = (\nu_1 \cdots \nu_n)^{-1} =
  1$. Again by Proposition \ref{prop:mult2boxes},~\eqref{eq:item4} has
  a non-trivial integer solution. By the same method we find that
  if~\eqref{eq:item4} has a non-trivial solution,
  then~\eqref{eq:item3} has a non-trivial solution.
\end{proof}

\begin{proof}[Proof of Theorem \ref{thm:multi-trans}]
  Assume that condition (i) of the theorem holds.  Let $q^{(j)}$ be
  the sequence of solutions to inequality \eqref{eq:15} and order
  these such that $j' > j$ implies that $F_+ (q^{(j')}) \geq F_+
  (q^{(j)})$, where
  \begin{equation*}
    F_+(x)=\left(\prod_{j=1}^n\max\{1,|x_j|\}\right)^{1/n}. 
  \end{equation*}
  Furthermore, since equality holds at most finitely often, and since
  there are infinitely many solutions, $F_+ (q^{(j)})$ tends to
  infinity as $j$ tends to infinity. Now for each $j \in \natnum$,
  define numbers
  \begin{equation*}
    \mu_j = F_+ (q^{(j)}), \quad \lambda_j = \mu_j^{-1 - \varepsilon}.
  \end{equation*}
  With these coefficients,~\eqref{eq:item3} has $q^{(j)}$ as a
  solution for each $j \in \natnum$. By Proposition
  \ref{prop:pre_transfer},~\eqref{eq:item4} has a non-trivial integer
  solution $p^{(j)}$ for each $j \in \natnum$. Let $j$ be fixed but
  arbitrary. For ease of notation, we drop the index in the following.

  From \eqref{eq:item4}, 
  \begin{equation*}
    F(\langle p a \rangle) \leq n \lambda, \quad \abs{p} \leq n
    \mu \lambda^{(1-n)/n}.
  \end{equation*}
  Let $\varepsilon' \in (0, \varepsilon)$ be fixed but arbitrary.  To
  prove that there is a solution to~\eqref{eq:16} of condition (ii),
  the $\lambda$ and $\mu$ are eliminated. The second inequality
  immediately implies that
  \begin{alignat}{2}
    \label{oldeq:8}
    \abs{p}^{-(1 + \varepsilon')/n} &\geq \left(n \mu
      \lambda^{(1-n)/n}\right)^{-(1+\varepsilon')/n} \notag\\
    &= \left(n^{-(1/n) - (\varepsilon'/n)-1} \mu^{-(1/n)
        -(\varepsilon'/n)} \lambda^{-(1/n^2) + (1/)n -
        (\varepsilon'/n^2) + (\varepsilon'/n) -1} \right)n \lambda.
\end{alignat}
  Hence, if we can prove that $\varepsilon'$ and $j$ may be chosen so
  that the term in the brackets is greater than or equal to $1$, this
  implies that we have a solution to inequality~\eqref{eq:16}. We
  insert the definition of $\lambda$ in \eqref{oldeq:8} to see that
  this is equivalent to the condition that
  \begin{equation*}
    n^{-(1/n)-(\varepsilon'/n)-1}\mu^{-(2/n)-(2\varepsilon'/n) +
      (1/n^2) + (\varepsilon'/n^2) + 1 + (\varepsilon/n^2) -
      (\varepsilon/n) + (\varepsilon \varepsilon'/n^2) - (\varepsilon 
      \varepsilon'/n) +\varepsilon} \ge 1.
  \end{equation*}
  Using the fact that $\mu$ is increasing and unbounded as a function
  of $j$ and on choosing $\varepsilon'$ sufficiently small, it is
  straightforward to verify the claim.
  
  To prove the first implication, we need to show that this produces
  infinitely many solutions to inequality~\eqref{eq:16}. But again, this
  follows from Proposition \ref{prop:pre_transfer}.  Indeed, assume
  that we have only solutions $p$ to \eqref{eq:item4} with bounded
  height. Then there are only solutions $q$ with bounded $F(q)$, and
  hence only finitely many solutions.

  The converse implication is proved analogously by interchanging the
  two inequalities in the above.
\end{proof}

The above proof can be generalised to other star bodies than the
multiplicative one. Indeed, we have only used the fact that solutions
to the relevant inequalities may be encoded in matrix
form~\eqref{oldeq:7}, parametrised by the numbers $\nu_i$.  This
parametrisation in turn depends only on exhibiting an appropriate
cover of the star body corresponding to $F$. In the case studied
above, the star body is covered by boxes with sides parallel to the
coordinate axes (see figure \ref{fig:multiplicative}). We now give an
example when this is not the case. For convenience, we restrict
ourselves to the planar case.

The `union jack' star body is obtained by taking a union of the
multiplicative star body and a rotation of it through $\pi/4$. The
associated distance function is given by 
\begin{equation*}
  F'(x,y) = \min\left\{\abs{xy}, \dfrac{\abs{x^2 -
        y^2}}{2}\right\}^{1/2}.
\end{equation*}
\begin{figure}[htbp]
  \centering
  \includegraphics[width=9cm]{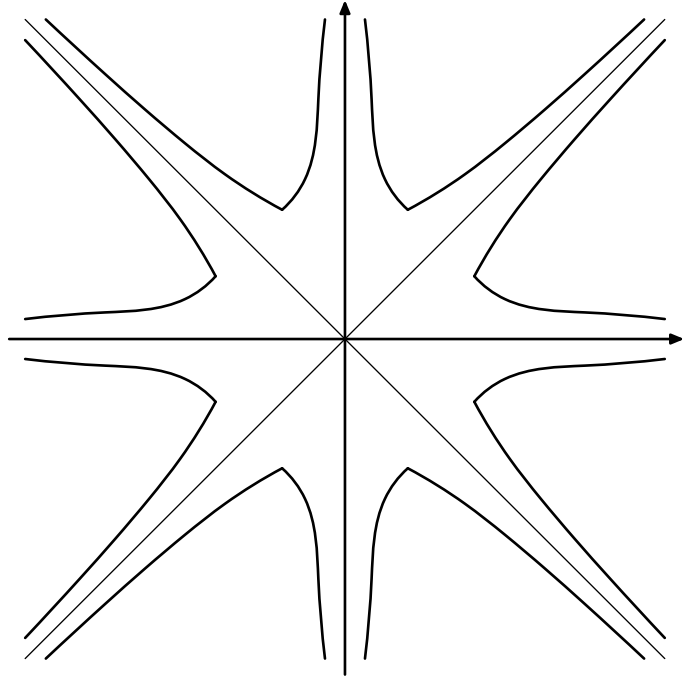}
  \caption{The `union jack' star body}
  \label{fig:union_jack}
\end{figure}
A plot of the star body is given in figure \ref{fig:union_jack}. For
this distance function, the above argument gives the following
transference principle.

\begin{thm}
  \label{thm:union-trans}
  Let $(x,y) \in (\mathbb{R}\setminus \rationals)^2$. The following
  two conditions are equivalent 
  \begin{enumerate}[(i)]
  \item For some $\varepsilon > 0$, there are infinitely many
    $q_1, q_2 \in \integer$ such that
    \begin{multline*}
      \abs{\langle q_1 x + q_2 y \rangle} \leq
      \min\Big\{\max\{\abs{q_1},1\} \max\{\abs{q_2},1\},\\
      \tfrac{\sqrt{2}}{2} \max\{\abs{q_1 + q_2},1\} \max\{\abs{q_1 - 
        q_2},1\}\Big\}^{-1-\varepsilon}.
    \end{multline*}
  \item For some $\varepsilon' > 0$, there are infinitely many
    $p \in \integer$ such that
    \begin{equation*}
      F'\left(\langle (p x, p y) \rangle \right) \leq
      \abs{p}^{-(1+\varepsilon')/2}. 
    \end{equation*}
  \end{enumerate}
\end{thm}

The proof is essentially the same, except that instead of matrices $A$
and $A^*$ from \eqref{oldeq:7}, we use four matrices
\begin{eqnarray*}
  \tilde{A} & = 
  \begin{pmatrix}
    \lambda^{-1}x & \nu_1 \mu^{-1} & 0 \\
    \lambda^{-1}y & 0 & \nu_2 \mu^{-1} \\
    \lambda^{-1} & 0 & 0
  \end{pmatrix},
  \quad
  \tilde{A}' &= 
  \begin{pmatrix}
    \lambda^{-1}x & \tfrac{\sqrt{2}}{2} \nu'_1 \mu^{-1} & -
    \tfrac{\sqrt{2}}{2} \nu_2' \mu^{-1} \\
    - \lambda^{-1}y & \tfrac{\sqrt{2}}{2} \nu_1' \mu^{-1} &
    \tfrac{\sqrt{2}}{2} \nu_2 \mu^{-1} \\
    \lambda^{-1} & 0 & 0
  \end{pmatrix}, \\
  \tilde{\tilde{A}} &=
  \begin{pmatrix}
    \lambda & - \mu \nu_1^{-1} x & -\mu \nu_2^{-1} y \\
    0 & \mu \nu_1^{-1} & 0 \\
    0 & 0 & \mu \nu_2^{-1} 
  \end{pmatrix},
  \quad \tilde{\tilde{A}}' &=
  \begin{pmatrix}
    \lambda & - \tfrac{\sqrt{2}}{2} \mu {\nu'_1}^{-1} (x-y) &
    \tfrac{\sqrt{2}}{2}  \mu {\nu'_2}^{-1}(x+y) \\
    0 & \tfrac{\sqrt{2}}{2} \mu {\nu'_1}^{-1} & \tfrac{\sqrt{2}}{2}
    \mu {\nu'_2}^{-1} \\ 
    0 & - \tfrac{\sqrt{2}}{2} \mu {\nu'_1}^{-1} &
    \tfrac{\sqrt{2}}{2}\mu {\nu'_2}^{-1} 
  \end{pmatrix},
\end{eqnarray*}
parametrised by $\nu_1, \nu_2, \nu_1', \nu_2' > 0$, where $\nu_1 \nu_2
= \nu_1'\nu_2' = 1$.  It is easily checked that these matrices encode
solutions to the inequalities corresponding to \eqref{eq:item3} and
\eqref{eq:item4} in the present setting and so provide an analogue of
Proposition \ref{prop:pre_transfer}, \emph{i.e.}, we now look at
integer solutions to
\begin{equation*}
  \min\left\{\infabs{\tilde{q} \tilde{A}}, \infabs{\tilde{q}
  \tilde{A}'} \right\} \leq 1,
\end{equation*}
and similarly to
\begin{equation*}
  \min\left\{\infabs{\tilde{p} \tilde{\tilde{A}}}, \infabs{\tilde{q}
    \tilde{\tilde{A}}'}\right\} \leq 1.
\end{equation*}
The proof of Theorem \ref{thm:union-trans} is now analogous to that of
Theorem \ref{thm:multi-trans}. We leave details to the reader.

Clearly, this approach can be generalised to other star bodies, where
the parametrisation depends on the distance function, and so the above
proof can be seen as a recipe for proving transference principles.
However, apart from the minor variations of standard known cases of
simultaneous (height), dual (inner product) and multiplicative
approximation (see \emph{e.g.} \cite{MR92h:11063, dyson47}), it does
not appear that there is a simple general transference principle for
arbitrary distance functions.

\section*{Acknowledgements}
\label{sec:acknowledgements}

We thank Victor Beresnevich for pointing out an error in an earlier
version of the paper and Yann Bugeaud for drawing our attention to
Dyson's paper~\cite{dyson47}. We are grateful for support from EPSRC
(grant no. GR/N02832/01) and INTAS (grant no.  001--429). SK is a
William Gordon Seggie Brown Fellow.

\end{document}